\documentclass[11pt]{article}
\usepackage{latexsym,amsfonts}
\usepackage{amssymb,amsmath}
\usepackage{enumitem}
\usepackage{graphics}
\usepackage[demo]{graphicx}
\usepackage{epsfig}
\usepackage{pstricks}
\usepackage[normal]{subfigure}
\usepackage[latin1]{inputenc}
\usepackage[english,francais]{babel}
\usepackage{relsize,exscale}
\usepackage{makeidx}
\usepackage{enumitem}
\usepackage{amsfonts,amssymb,amsmath}
\usepackage{graphicx}
\usepackage{color}
\usepackage{multirow}
\usepackage{mathrsfs}
\usepackage[normalem]{ulem}
\usepackage{cancel}
\usepackage{bbm}
\newenvironment{prooff}{{\it Proof :}}{\hfill\rule{2mm}{2mm}\vskip3mm\par}
\newtheorem{theorem}{Theorem}[section]
\newtheorem{lemma}[theorem]{Lemma}

\newtheorem{e-definition}[theorem]{Definition\rm}
\newtheorem{remark}{\it Remark\/}

\usepackage{color}
\definecolor{dred}{rgb}{0.92,0,0}
\definecolor{dgreen}{rgb}{0,0.92,0}
\definecolor{dblue}{rgb}{0,0,0.92}
\definecolor{dyellow}{rgb}{0.95,0.95,0}

\newcommand{\R}{\mathbb{R}}
\newcommand{\N}{\mathbb{N}}
\def\D{\displaystyle}
\newcommand{\ds}{\displaystyle}
\newcommand{\hs}{\hspace{0.1cm}}

\newcommand{\sa}{\\ [0.2cm]}
\usepackage{multirow}
\graphicspath{
{./Figures/}
{Figures/}
{./}
}
\title{A new second order Taylor-like theorem \\ with an optimized reduced remainder}
\author{Jo\"el Chaskalovic \thanks{D'Alembert,
Sorbonne University, Paris, France, (Email:  jch1826@gmail.com)}
\qquad
Franck Assous
\thanks{
Ariel University, Ariel, Israel, (Email: \emph{(corresp.)} assous@ariel.ac.il)}
\qquad
Hessam Jamshidipour
\thanks{
D'Alembert,
Sorbonne University, Paris, France, (Email:  hessam.jamshidipour@gmail.com)}
}
\date{}
\begin{document}
\maketitle
\selectlanguage{english}
\begin{abstract}
\noindent In this paper, we derive a variant of the Taylor theorem to obtain a new minimized remainder. For a given function $f$ defined on the interval $[a,b]$, this formula is derived by introducing a linear combination of $f'$ computed at $n+1$ equally spaced points in $[a,b]$, together with  $f''(a)$ and $f''(b)$. We then consider two classical applications of this Taylor-like expansion: the interpolation error and the numerical quadrature formula. We show that using this approach improves both the Lagrange $P_2$- interpolation error estimate and the error bound of the Simpson rule in numerical integration.
\end{abstract}
\noindent {\em keywords}: Taylor's theorem, Lagrange interpolation, interpolation error, Simpson rule, quadrature error.
\section{Introduction}\label{intro}

\noindent Even today, improving the accuracy of approximation remains a challenging problem in numerical analysis. Here, we are concerned  with the difficulty of accurately determining the error estimate in numerical methods. This article is part of a series of articles in which this topic is addressed. More precisely, we derive here a variant of Taylor's theorem to obtain a new minimized remainder that we apply to  interpolation error and numerical quadrature formula.\\

\noindent From a mathematical point of view,  the origin of such problems already appears in Rolle's theorem and in Lagrange and Taylor's theorems, see for instance \cite{Atki88}. Indeed, there exists an unknown point involved in the remainder of Taylor's expansion, that leads to some ``uncertainty''.\\

\noindent Consequently,  most error estimates focus generally on the {\em asymptotic behavior} of the error. For instance in finite element approximation, {\em a priori} error estimates consider the asymptotic behavior of the difference between the exact and the approximate solution, as the mesh size h tends to zero.\\

\noindent  In this context, several approaches have been proposed to determine a way to improve the accuracy of approximation. For example, within the framework of numerical integration, we refer the reader to \cite{Barnett_Dragomir}, \cite{Cerone} or \cite{Dragomir_Sofo}, and references therein. From another point of view, due to the lack of information, heuristic methods were considered, basically based on a probabilistic approach, see for instance \cite{Abdulle}, \cite{AsCh2014}, \cite{Hennig}, \cite{Oates} or \cite{ArXiv_JCH}, \cite{CMAM2} and \cite{ChAs20}. This allows to compare different numerical methods, and more precisely finite element, for a given fixed mesh size, see \cite{MMA2021}.\\

\noindent Nevertheless, Taylor's theorem introducing an unknown point in its remainder, this makes it very difficult (often impossible) to compute the interpolation error, and consequently the approximation error of a given numerical method. Therefore, the possibility to accurately estimate the upper bounds of the error remains an important issue. In this article, we study the values of the numerical constants involved in such estimates, trying to reduce them as small as possible.\\

\noindent In this framework, we proposed in \cite{ChAs2023} a refined first-order expansion formula in $\R^n$,  to get an reduced remainder, compared to the one obtained by usual Taylor's formula. Then, we investigate the related properties in the interpolation error estimates and in Lagrange finite element error estimates. In the context of quadrature rules applications, such a problem was considered in the past years, and is often referred as the perturbed (or corrected) quadrature rules, see for instance \cite{Cerone} or \cite{Dragomir_Sofo}.  In other examples, the authors obtained in \cite{Chen01}, \cite{DrWa97} or \cite{MaPU00}, the trapezoid inequality by the difference of sup and inf bound of the first derivative.\\

\noindent In this paper, we are concerned by a second order Taylor-like theorem, leading to an optimized reduced remainder. Applications we have in mind are the interpolation error based on a second-order polynomial and the Simpson quadrature rule (see for instance \cite{BuFa11}). Concerning the Simpson inequality, we also refer the reader to \cite{Liu02}.\\

\noindent The main difficulty addressed in this article concerns the development and use of a new Taylor-type formula, with the smallest possible remainder. Another important aspect is the application of this new formula to interpolation error and numerical quadrature formulas. The paper is organized as follows. In Section \ref{Section-newsecond}, we present the main result, which treats on the improved second order Taylor-like formula. In Section \ref{C}, we consider two classical applications of the Taylor expansion: the interpolation error is investigated in subsection \ref{Interpolation_Error}, whereas the Simpson's quadrature rule is studied subsection \ref{Quadrature_Error}. In both cases, we derive new results on error estimate. Concluding remarks follow.

\section{A new second order expansion formula}\label{Section-newsecond}
\noindent To begin with, let us recall the well known second order Taylor's formula \cite{Taylor}. We consider $(a,b) \in \mathbb{R}^{2}$, $a < b$, and a function $f \in \mathcal{C}^{3}([a,b])$. Then, following Taylor's theorem in one real variable, there exist two real constants $m_{3}$ and $M_{3}$ such that, for all $x \in [a,b]$,
\begin{equation}\label{m3M3}
\D -\infty < m_3=\inf_{a\leq x\leq b}f'''(x) \mbox{ and } \D M_3=\sup_{a\leq x\leq b}f'''(x) < +\infty\,,
\end{equation}
and we have
\begin{equation}\label{Second_Order_Taylor}
f(b) = f(a) + (b-a)f'(a) +  \frac{(b-a)^2}{2} f''(a) +(b-a)^2\epsilon_{a,2}(b),
\end{equation}
where
\[ \lim_{b \to a} \epsilon_{a,2}(b) = 0, \]
and
\begin{equation}\label{Epsilon_b_bounded}
\frac{(b-a)}{6}m_{3} \leqslant \D \epsilon_{a,2}(b) \leqslant \frac{(b-a)}{6}M_{3}.
\end{equation}

\noindent In the same spirit we proposed in \cite{arXiv_First_Order} for the first order case, our aim is now to derive a new second order Taylor-like formula that gives us a minimized remainder. To that aim, let us first recall the main result obtained in \cite{arXiv_First_Order} for the first order case. Given two reals $a,b \in \mathbb{R}$, $a < b$ and an integer $n \in \mathbb{N}^{*}$, we proved the following result
\begin{theorem}\label{theorem_1}
Let $f$ be a real mapping define on $[a,b]$ which belongs to $\mathcal{C}^{2}([a,b])$, such that: $\forall x \in [a,b], -\infty < m_2 \leqslant f''(x) \leqslant M_2 < +\infty$. \sa
Then we have the following first order expansion:
\begin{equation}\label{Generalized_Taylor1}
f(b) = f(a) + (b-a)\left(\frac{f'(b) + f'(a)}{2n} + \frac{1}{n}\sum \limits_{k=1}^{n-1} f'\left(a + k\frac{(b-a)}{n}\right)\right) + (b-a)\epsilon_{a,n+1}^{(1)}(b),
\end{equation}
where :
$$
\D|\epsilon_{a,n+1}^{(1)}(b)| \leqslant \frac{(b-a)}{8n}(M_2-m_2).
$$
\end{theorem}
Moreover, this result is optimal in the sense that the weights in (\ref{Generalized_Taylor1}) involved in the linear combination of $f'$ at the equally spaced points $\D a + k\frac{(b-a)}{n}$ guarantee the remainder $\epsilon_{a,n+1}^{(1)}(b)$ to be minimum.\\

\noindent In order to prove the main theorem for the second order case considered in this paper, we will need the following lemma proved in \cite{arXiv_First_Order}:
\begin{lemma} \label{formula_1}
Let $u$ be a continuous function on $[a,b]$, and, for $n\in\N^*$,  let $(a_{k})_{0 \leqslant k \leqslant n}$ be a finite sequence of real numbers. We have the following formula:
$$
 \sum \limits_{k=0}^{n-1} \int_{k}^{n}{a_{k} u(t) dt} = \sum \limits_{k=0}^{n-1} \int_{k}^{k+1}{S_{k} u(t) dt}, \quad \mbox {with } S_{k} = \sum \limits_{j=0}^{k} a_{j}.
$$
\end{lemma}
From now on, we assume that $n \in \mathbb{N}^{*}$. To obtain a second order Taylor-like formula, we first consider the following generalization of (\ref{Generalized_Taylor1}) involving the ``reminder'' $\epsilon_{a,n+1}^{(2)}(b)$:
\begin{eqnarray}
&&\D f(b) = f(a) + (b-a)\left(\frac{f'(b) + f'(a)}{2n} + \frac{1}{n}\sum \limits_{k=1}^{n-1} f'\left(a + k\frac{(b-a)}{n}\right) \right) \nonumber\\
&&\hspace*{2.cm}+ (b-a)^2\ \sum \limits_{k=0}^{n} \omega_{k}(n)f''\left(a + k\frac{(b-a)}{n}\right) + (b-a)^2\epsilon_{a,n+1}^{(2)}(b),\label{aftL22}
\end{eqnarray}
that we rewrite for simplicity as
$$
\D f(b) = f(a) + (b-a)\Lambda_{n}^{(1)}(a,b)
+ (b-a)^2\Lambda_{n}^{(2)}(a,b) + (b-a)^2\epsilon_{a,n+1}^{(2)}(b),
$$
with
\begin{equation}
\label{eq-Lambda1}
\D \Lambda_{n}^{(1)}(a,b) = \frac{f'(b) + f'(a)}{2n} + \frac{1}{n}\sum \limits_{k=1}^{n-1} f'\left(a + k\frac{(b-a)}{n}\right),
\end{equation}
and $\Lambda_{n}^{(2)}(a,b)$ defined by determining
\begin{equation}\label{Lambda2}
\Lambda_{n}^{(2)}(a,b) = \sum \limits_{k=0}^{n} \omega_{k}(n)f''\left(a + k\frac{(b-a)}{n}\right)\,,
\end{equation}
such that (\ref{aftL22}) holds, and so, $\epsilon_{a,n+1}^{(2)}(b)$ goes to 0 when $b \rightarrow a$.
Our aim is now to determine the sequence of real weights $(\omega_{k}(n))_{0 \leq k \leq n}$ that minimizes  the remainder $\epsilon_{a,n+1}^{(2)}(b)$. this result is stated in the following theorem:
\begin{theorem}\label{theorem_2}
Let $f$ be a real mapping defined on $[a,b]$ which belongs to $\mathcal{C}^{3}([a,b])$, such that: $\forall x \in [a,b], -\infty < m_3 \leqslant f'''(x) \leqslant M_3 < +\infty$. If the weights $(\omega_k(n))_{0 \leq k \leq n}$ satisfy $\D \sum \limits_{k=0}^{n} \omega_{k}(n) = 0$,
then we have:
\begin{equation}\label{Poids}
\omega_{0}(n) = -\omega_{n}(n) = \frac{3}{32n^2} \hs \mbox{ and } \hs \omega_{k}(n) = 0, \hs \forall 0 < k < n,
\end{equation}
and the following second order expansion hold:
\begin{equation}\label{Generalized_Taylor}
\D f(b) = f(a) + (b-a)\Lambda_{n}^{(1)}(a,b) + (b-a)^2\Lambda_{n}^{(2)}(a,b) + (b-a)^2\epsilon_{a,n+1}^{(2)}(b),
\end{equation}
where $\D \Lambda_{n}^{(1)}(a,b)$ is given by (\ref{eq-Lambda1}) and $\D \Lambda_{n}^{(2)}(a,b)$ is expressed as
\begin{eqnarray}
\label{L2}
\D \Lambda_{n}^{(2)}(a,b) & = & -\frac{3}{32n^2}\biggl(f''(b)-f''(a)\biggr)\,.
\end{eqnarray}
Moreover, this result is optimal  since the weights introduced in (\ref{Poids}) guarantee that the remainder $\epsilon_{a,n+1}^{(2)}(b)$ is minimum, and satisfies:
\begin{equation}\label{Control_Epsilon2-base}
\D \frac{(b-a)}{96n^2}(2m_3-M_3) \leq \epsilon_{a,n+1}^{(2)}(b) \leq \frac{(b-a)}{96n^2}(2M_3-m_3).
\end{equation}
Consequently, $\displaystyle \lim_{b \rightarrow a} \epsilon_{a,n+1}^{(2)}(b)=0$.
\end{theorem}
\begin{prooff}
Let us observe first that, using (\ref{Generalized_Taylor1}), $\epsilon_{a,n+1}^{(2)}(b)$ can be written as
\begin{equation}\label{eps2exp}
\D\epsilon_{a,n+1}^{(2)}(b) = \frac{\epsilon_{a,n+1}^{(1)}(b)}{b-a} - \Lambda_{n}^{(2)}(a,b)\,.
\end{equation}
In \cite{arXiv_First_Order}, it is proved (see expression (16) together with (28)) that the remainder $\epsilon_{a,n+1}^{(1)}(b)$ of the expansion (\ref{Generalized_Taylor1})
 can be expressed by
\begin{equation}\label{F1}
\D\epsilon_{a,n+1}^{(1)}(b) = \sum \limits_{k=0}^{n-1} \int_{\frac{k}{n}}^{\frac{k+1}{n}}{\biggl(\frac{1}{2n}+\frac{k}{n} - t\biggr)\phi'(t)dt},
\end{equation}
where $\phi'$ is the derivative of the function $\phi$ defined by:
$$
\begin{array}{r c c c l}
    \phi & : & [0,1] & \longrightarrow & \mathbb{R} \\
         &   &     t & \longmapsto & f'(a + t(b-a)).
\end{array}
$$
Now, we perform an integration by parts of the integral involved in (\ref{F1}) and we get
\begin{eqnarray*}
\D \int_{\frac{k}{n}}^{\frac{k+1}{n}}\hspace{-0.1cm}\biggl(\frac{1}{2n}+\frac{k}{n} - t\biggr)\phi'(t)dt \hspace{-0.2cm}& = &
\hspace{-0.2cm}\frac{k(k+1)}{2n^2}\biggl[\phi'\biggl(\frac{k+1}{n}\biggr)-\phi'\biggl(\frac{k}{n}\biggr)\biggr]
- \int_{\frac{k}{n}}^{\frac{k+1}{n}}{\biggl[\biggl(\frac{1}{2n}+\frac{k}{n}\biggr)t - \frac{t^2}{2}\biggr]\phi''(t)dt}, \nonumber\\[0.2cm]
& = & \int_{\frac{k}{n}}^{\frac{k+1}{n}}{\biggl[\frac{k(k+1)}{2n^2}-\biggl(\frac{1}{2n}+\frac{k}{n}\biggr)t + \frac{t^2}{2}\biggr]\phi''(t)dt}.
\end{eqnarray*}
Using this expression, $\epsilon_{a,n+1}^{(1)}(b)$ can be rewritten as:
$$
\D\epsilon_{a,n+1}^{(1)}(b) = \sum \limits_{k=0}^{n-1} \int_{\frac{k}{n}}^{\frac{k+1}{n}}{\biggl[\frac{t^2}{2}-\biggl(\frac{1}{2n}+\frac{k}{n}\biggr)t
+\frac{k(k+1)}{2n^2}\biggr]\phi''(t)dt}.
$$
Consequently, from (\ref{eps2exp}), $\epsilon_{a,n+1}^{(2)}(b)$ can be written
\begin{equation}\label{F3}
\D\epsilon_{a,n+1}^{(2)}(b) = \frac{1}{b-a}\sum \limits_{k=0}^{n-1} \int_{\frac{k}{n}}^{\frac{k+1}{n}}{\biggl[\frac{t^2}{2}-\biggl(\frac{1}{2n}+\frac{k}{n}\biggr)t
+\frac{k(k+1)}{2n^2}\biggr]\phi''(t)dt} - \frac{1}{b-a}\sum \limits_{k=0}^{n}\omega_k(n)\phi'\biggl(\frac{k}{n}\biggr).
\end{equation}
Now,  using that $\D\phi'\biggl(\frac{k}{n}\biggr)= \phi'(1)-\int_{\frac{k}{n}}^{1}\phi''(t)dt$ together with  Lemma \ref{formula_1}, and a change of variable $t'=n t$, the last sum in (\ref{F3}) can be expressed as:
$$
\D\sum \limits_{k=0}^{n}\omega_k(n)\phi'\biggl(\frac{k}{n}\biggr) = \sum \limits_{k=0}^{n}\omega_k(n)\phi'(1)
- \sum \limits_{k=0}^{n-1}\int_{\frac{k}{n}}^{\frac{k+1}{n}}S_k(n)\phi''(t)dt,
$$
where $\D S_k(n)=\sum_{j=0}^{k}\omega_j(n),$ for all $0 \leq k \leq n-1$.\sa
Consequently, expression (\ref{F3}) becomes:
\begin{equation}\label{F5}
\epsilon_{a,n+1}^{(2)}(b) = \frac{1}{b-a}\sum \limits_{k=0}^{n-1} \int_{\frac{k}{n}}^{\frac{k+1}{n}}{\biggl[\frac{t^2}{2}-\biggl(\frac{1}{2n}+\frac{k}{n}\biggr)t
+\frac{k(k+1)}{2n^2}+S_k(n)\biggr]\phi''(t)dt} -  \frac{1}{b-a}\phi'(1)\sum \limits_{k=0}^{n}\omega_k(n).
\end{equation}
Let us assume for simplicity (see Remark \ref{Relation_Poids} below) that:
\begin{equation}\label{poids_norme}
\D \sum \limits_{k=0}^{n} \omega_{k}(n) = 0.
\end{equation}
Then, setting
$$
\D\lambda=S_k(n)+\frac{k(k+1)}{2n^2},
$$
the expression (\ref{F5}) of $\epsilon_{a,n+1}^{(2)}(b) $ becomes:
$$
\epsilon_{a,n+1}^{(2)}(b) =\frac{1}{b-a} \sum \limits_{k=0}^{n-1} \int_{\frac{k}{n}}^{\frac{k+1}{n}}{\biggl[\frac{t^2}{2}-\biggl(\frac{1}{2n}+\frac{k}{n}\biggr)t+
\lambda\biggr]\phi''(t)dt}.
$$
Substituting $\D t=\frac{k+s}{n}$ in this integral, and setting $P_{\bar{\lambda}}(s) = s^2 - s+ \bar{\lambda}$ with
\begin{equation}\label{lambdabar}
\bar{\lambda} \equiv 2n^2\lambda -k(k+1) = 2n^2 S_k(n)\,,
\end{equation}
we get:
\begin{equation}\label{F7}
\epsilon_{a,n+1}^{(2)}(b) = \frac{1}{2(b-a)n^3}\sum \limits_{k=0}^{n-1} \int_{0}^{1}P_{\bar{\lambda}}(s)\phi''\biggl(\frac{s+k}{n}\biggr)ds.
\end{equation}
Now, assuming that the discriminant $\Delta=1-4\bar{\lambda}$ of $P_{\bar{\lambda}}$ is strictly positive, it exists $(t_1,t_2)\in\R^2, t_1 < t_2,$ such that: $P_{\bar{\lambda}}(t_1) = P_{\bar{\lambda}}(t_2)=0$.\sa
In the following, our aim is to derive an estimate of $\epsilon_{a,n+1}^{(2)}(b)$. For this purpose, we split  the integral above depending on the roots of $P_{\bar{\lambda}}$. We get:
$$
\D \int_{0}^{1}P_{\bar{\lambda}}(t)\phi''\biggl(\frac{t+k}{n}\biggr)dt = \int_{0}^{t_1}P_{\bar{\lambda}}(t)\phi''\biggl(\frac{t+k}{n}\biggr)dt
+ \int_{t_1}^{t_2}P_{\bar{\lambda}}(t)\phi''\biggl(\frac{t+k}{n}\biggr)dt + \int_{t_2}^{1}P_{\bar{\lambda}}(t)\phi''\biggl(\frac{t+k}{n}\biggr)dt.
$$
Keeping in mind that, for all $x \in [0,1]$ and for all $t\in [a,b]$,
$$
\phi''(x) = (b-a)^2 f'''(a + x(b-a)) \hs \mbox{ and }  \,\,m_3 \leqslant f'''(t) \leqslant M_3,
$$
then, $P_{\bar{\lambda}}(t)$ keeps a constant sign on each of the three above integrals and we have:
\begin{eqnarray*}
\D m_3(b-a)^2\int_{0}^{t_1}P_{\bar{\lambda}}(t)dt & \leq & \int_{0}^{t_1}P_{\bar{\lambda}}(t)\phi''\biggl(\frac{t+k}{n}\biggr)dt \leq M_3(b-a)^2\int_{0}^{t_1}P_{\bar{\lambda}}(t)dt, \label{I.1}\\[0.2cm]
\D M_3(b-a)^2\int_{t_1}^{t_2}P_{\bar{\lambda}}(t)dt & \leq & \int_{t_1}^{t_2}P_{\bar{\lambda}}(t)\phi''\biggl(\frac{t+k}{n}\biggr)dt \leq m_3(b-a)^2\int_{t_1}^{t_2}P_{\bar{\lambda}}(t)dt, \label{I.2} \\[0.2cm]
\D m_3(b-a)^2\int_{t_2}^{1}P_{\bar{\lambda}}(t)dt & \leq & \int_{t_2}^{1}P_{\bar{\lambda}}(t)\phi''\biggl(\frac{t+k}{n}\biggr)dt \leq M_3(b-a)^2\int_{t_2}^{1}P_{\bar{\lambda}}(t)dt,\label{I.3}
\end{eqnarray*}
that yields the two following inequalities
\begin{eqnarray}
\D \int_{0}^{1}P_{\bar{\lambda}}(t)\phi''\biggl(\frac{t+k}{n}\biggr)dt \hspace{-0.2cm}& \leq & \hspace{-0.2cm}m_3(b-a)^2\int_{t_1}^{t_2}\hspace{-0.2cm}P_{\bar{\lambda}}(t)dt +
M_3(b-a)^2\biggl[\int_{0}^{t_1}\hspace{-0.2cm}P_{\bar{\lambda}}(t)dt+\int_{t_2}^{1}\hspace{-0.2cm}P_{\bar{\lambda}}(t)dt\biggr], \label{I.4}\\[0.2cm]
\D \int_{0}^{1}P_{\bar{\lambda}}(t)\phi''\biggl(\frac{t+k}{n}\biggr)dt \hspace{-0.2cm}& \geq &\hspace{-0.2cm} M_3(b-a)^2\int_{t_1}^{t_2}\hspace{-0.2cm}P_{\bar{\lambda}}(t)dt +
m_3(b-a)^2\biggl[\int_{0}^{t_1}\hspace{-0.2cm}P_{\bar{\lambda}}(t)dt+\int_{t_2}^{1}\hspace{-0.2cm}P_{\bar{\lambda}}(t)dt\biggr].\label{I.5}
\end{eqnarray}
We have now to deal with these inequalities. Since they have the same structure, we will consider only the first one, the second one can be treated in the same way. \sa
Dividing by $(b-a)^2$ and computing the integrals, $P_{\bar{\lambda}}(t)$ being a second-degree polynomial function, we easily get:
\begin{eqnarray*}
\D \hspace{-0.4cm}\frac{\D\int_{0}^{1}\hspace{-0.2cm}P_{\bar{\lambda}}(t)\phi''\biggl(\frac{t+k}{n}\biggr)dt}{(b-a)^2} \hspace{-0.2cm}& \leq &
\hspace{-0.3cm}(M_3-m_3)(t_1-t_2)\biggl[\frac{t_1^2+t_1t_2+t_2^2}{3}-\frac{t_1+t_2}{2}+\bar{\lambda}\biggr]+M_3\biggl(\bar{\lambda}-\frac{1}{6}\biggr).
\end{eqnarray*}
Using that $t_i,(i=1,2)$, are the roots of the polynomial $P_{\bar{\lambda}}(.)$, we have, for $i=1,2$,  $t_i^2=t_i-\bar{\lambda}$ and this inequality becomes:
$$
\frac{\D\int_{0}^{1}\hspace{-0.2cm}P_{\bar{\lambda}}(t)\phi''\biggl(\frac{t+k}{n}\biggr)dt}{(b-a)^2} \leq
(M_3-m_3)(t_1-t_2)\biggl[\frac{t_1+t_2+t_1t_2-2\bar{\lambda}}{3}-\frac{t_1+t_2}{2}+\bar{\lambda}\biggr]+M_3\biggl(\bar{\lambda}-\frac{1}{6}\biggr).
$$
Since $t_1, t_2$ are the roots of the second-degree polynomial $P_{\bar{\lambda}}(.)$,
$$
t_1+t_2 = 1, \hs t_1t_2 = \bar{\lambda} \hs \mbox{ and } \hs t_1-t_2= -\sqrt{1-4\bar{\lambda}}\,,
$$
that leads to:
$$
\frac{\D\int_{0}^{1}\hspace{-0.2cm}P_{\bar{\lambda}}(t)\phi''\biggl(\frac{t+k}{n}\biggr)dt}{(b-a)^2} \leq
\frac{(M_3-m_3)}{6}(1-4\bar{\lambda})^{3/2} + M_3\biggl(\bar{\lambda}-\frac{1}{6}\biggr).
$$
Finally, due to the symmetry between $m_3$ and $M_3$ in the right-hand sides of (\ref{I.4}) and (\ref{I.5}), we can write
\begin{equation} \label{I.10}
\D \varphi_1(\bar{\lambda}) \leq \frac{\D\int_{0}^{1}\hspace{-0.2cm}P_{\bar{\lambda}}(t)\phi''\biggl(\frac{t+k}{n}\biggr)dt}{(b-a)^2} \leq \varphi_2(\bar{\lambda}),
\end{equation}
where $\varphi_i(\bar{\lambda}), (i=1,2)$ are defined by:
\begin{eqnarray}
\varphi_1(\bar{\lambda}) & = & \frac{(m_3-M_3)}{6}(1-4\bar{\lambda})^{3/2} + m_3\biggl(\!\bar{\lambda}-\frac{1}{6}\biggr), \nonumber\\[0.2cm]
\varphi_2(\bar{\lambda}) & = & \frac{(M_3-m_3)}{6}(1-4\bar{\lambda})^{3/2} + M_3\biggl(\!\bar{\lambda}-\frac{1}{6}\biggr).\nonumber
\end{eqnarray}
We conclude the proof by determining the $\bar{\lambda}$ that minimizes the distance between $\varphi_1(\bar{\lambda})$ and $\varphi_2(\bar{\lambda})$. Let us define
$$\varphi(\bar{\lambda}) = \varphi_2(\bar{\lambda})-\varphi_1(\bar{\lambda}) = (M_3-m_3)\biggl[\frac{(1-4\bar{\lambda})^{3/2}}{3} + \bar{\lambda}-\frac{1}{6}\biggr]\,,
$$
which satisfies $\varphi'(\bar{\lambda})=0$ for $\bar{\lambda} = \D\frac{3}{16}$, that is,  the minimum of $\varphi(\bar{\lambda})$. This also shows {\em a posteriori} that $1-4\bar{\lambda}>0$, i.e. the discriminant of $P_{\bar{\lambda}}(t)$ is positive. For this value of $\bar{\lambda}$, the inequalities (\ref{I.10}) are written as:
$$
\D \frac{2m_3-M_3}{48} \leq \frac{\D\int_{0}^{1}\hspace{-0.2cm}P_{\bar{\lambda}}(t)\phi''\biggl(\frac{t+k}{n}\biggr)dt}{(b-a)^2} \leq \frac{2M_3-m_3}{48},
$$
and summing then over $k$, we obtain for $ \epsilon_{a,n+1}^{(2)}(b)$ in (\ref{F7}):
$$
\D \frac{(b-a)}{96n^2}(2m_3-M_3) \leq \epsilon_{a,n+1}^{(2)}(b) \leq \frac{(b-a)}{96n^2}(2M_3-m_3).
$$
Moreover, we can also get the weights $\omega_k(n), (k=0,n),$ involved in $\Lambda_{n}^{(2)}(a,b)$ (cf. (\ref{Lambda2})). Indeed, substituting
$\D\bar{\lambda}=\frac{3}{16}$ in the expression $\D\bar{\lambda}=2n^2 S_k(n)$, from (\ref{lambdabar}), we get
$$
\D S_k(n) = \sum_{j=0}^{k}\omega_j(n) = \frac{3}{32n^2}, \forall k=0,\dots,n-1.
$$
Hence, for $k=0$, $\D\omega_0(n)=\frac{3}{32n^2}$ whereas $\omega_k(n)=0$ for all $0 \leq  k \leq n-1$. Finally, we determine the last weight $w_n(n)$ by using the assumption (\ref{poids_norme}), that leads to
$$
\D\omega_n(n) = -\omega_0(n) = -\frac{3}{32n^2}.
$$
and $\Lambda_{n}^{(2)}(a,b)$ introduced in (\ref{Lambda2}) satisfies (\ref{L2}).
\end{prooff}
\begin{remark}\label{Relation_Poids}
Condition (\ref{poids_norme}) on the weights $\omega_k(n)$ in Theorem \ref{theorem_2} is a kind of {\em closure condition} but is not a restrictive one. Indeed, without (\ref{poids_norme}), we will replace in (\ref{F7}) $\epsilon_{a,n+1}^{(2)}(b)$ by
$$
\epsilon_{a,n+1}^{(2)}(b) = \frac{1}{2(b-a)n^3}\sum \limits_{k=0}^{n-1} \int_{0}^{1}P_{\bar{\lambda}}(t)\phi''\biggl(\frac{t+k}{n}\biggr)dt - \frac{1}{b-a}\phi'(1)\sum \limits_{k=0}^{n-1}\omega_k(n)\,,
$$
and consequently,  the corresponding weights $\omega_k(n)$ would be written as:
$$
\omega_{0}(n) = \frac{3}{32n^2} \hs \mbox{ and } \hs \omega_{k}(n) = 0, \hs \forall\, 1 \leq k \leq n.
$$
Then, following the same steps, we will get the following estimates for $ \epsilon_{a,n+1}^{(2)}(b)$:
$$
\D  \frac{(2m_3-M_3)(b-a)}{96n^2}-\frac{3}{32n^2}\frac{\phi'(1)}{(b-a)} \leq \epsilon_{a,n+1}^{(2)}(b) \leq  \frac{(2M_3-m_3)(b-a)}{96n^2}-\frac{3}{32n^2}\frac{\phi'(1)}{(b-a)},
$$
or also
$$
\D \frac{(2m_3-M_3)(b-a)}{96n^2}-\frac{3M_2}{32 n^2} \leq \epsilon_{a,n+1}^{(2)}(b) \leq \frac{(2M_3-m_3)(b-a)}{96n^2}-\frac{3m_2}{32 n^2}.
$$
Furthermore, still without the condition (\ref{poids_norme}) on the weights, the second order Taylor's-like formula (\ref{Generalized_Taylor}) would be expressed as
\begin{equation}\label{Generalized_Taylor2_V2}
\D f(b) = f(a) + (b-a)\Lambda_{n}^{(1)}(a,b) +\frac{3}{32n^2} (b-a)^2f''(a) + (b-a)^2\epsilon_{a,n+1}^{(2)}(b),
\end{equation}
the coefficient before $f''(b)$ vanishes.
\end{remark}
Let us compare now the remainder (\ref{Control_Epsilon2-base}) of the new formula (\ref{Generalized_Taylor})-(\ref{L2}) with the reminder (\ref{Epsilon_b_bounded}) of the classical formula (\ref{Second_Order_Taylor}). As one can see, the former remainder is significantly smaller than the latter one: indeed, we have to compare $1/6$ with $1/32 n^2$ whose ratio is equal $3/16 n^2$. The worst case of this ratio corresponds to $n=1$ where the new remainder is approximatively 5 times smaller than the one obtained by the classical  formula.
\section{Applications to the approximation error}\label{C}

\noindent In this section, we consider two classical applications of the Taylor expansion: the Lagrange polynomial  interpolation and the numerical quadrature. In both cases, we will derive new formula of interpolation and quadrature, obtained by using the refined second-order expansion formula (\ref{Generalized_Taylor})-(\ref{L2}). Then, we will compare the errors obtained when using the standard Taylor expansion and our generalized approach. We begin with the interpolation error.

\subsection{The interpolation error}\label{Interpolation_Error}
\noindent Consider first the generalized Taylor-like expansion (\ref{Generalized_Taylor})-(\ref{L2}) for $n=2$ In this case, for any function $f$ which belongs to $C^3([a,b])$, this formula is expressed as
\begin{equation}\label{Generalized_Taylor-bis}
\D f(b) = f(a) + (b-a) \frac{f'(b) +2 f'\bigg(\!\ds\frac{a + b}{2}\!\bigg)+f'(a)}{4}
-\frac{3 (b-a)^2}{128}\left(f''(b)-f''(a)\right)
+ (b-a)^2\epsilon_{a,3}^{(2)}(b),
\end{equation}
where the remainder $\epsilon_{a,3}(b)$ satisfies
\begin{equation}\label{Control_Epsilon2}
\D \frac{(b-a)}{384}(2m_3-M_3) \leq \epsilon_{a,3}^{(2)}(b) \leq \frac{(b-a)}{384}(2M_3-m_3).
\end{equation}
\noindent In order to derive a first application of this formula, let us consider the case of the $P_2$-Lagrange interpolation (see for instance \cite{BuFa11}, \cite{Atki88}), where a given function $f$ is interpolated on $[a,b]$ by a polynomial $\Pi_{[a,b]}(f)$ of degree less than or equal to two. Hence, we can write:
\begin{equation}\label{P2_Interpolate}
\forall x \in [a,b]: \Pi_{[a,b]}(f)(x) =\frac{(x-c)(x-b)}{(a-c)(a-b)} \hs f(a) +  \frac{(x-a)(x-c)}{(b-a)(b-c)} \hs f(b)  +  \frac{(x-a)(x-b)}{(c-a)(c-b)} \hs f(c)\,,
\end{equation}
where $c=\ds\frac{a+b}{2}$ denotes the midpoint of $[a,b]$.\sa
As it is well known, one has by construction $\Pi_{[a,b]}(f)(a) = f(a)$, $\Pi_{[a,b]}(f)(b) = f(b)$ and $\Pi_{[a,b]}(f)(c) = f(c)$.
\sa
Let us investigate the consequences of formula (\ref{Generalized_Taylor-bis}) when it is used to evaluate the interpolation error $e(.)$ defined by
$$
\forall x \in [a,b]: e(x) = \Pi_{[a,b]}(f)(x) - f(x)\,,
$$
and compare it with the one obtained when using the classical second order Taylor's formula.\sa
The classical result \cite{Crouzeix_Mignot}, \cite{RaTho82} concerning the $P_2-$Lagrange interpolation error claims that, for any function $f$ that belongs to $C^3([a,b])$, we have:
$$
\D|e(x)| \leq \frac{(x-a)(b-x)|x-c|}{6}\sup_{a\leq x\leq b}|f'''(x)|\,,
$$
and using that
$
\D \sup_{a\leq x\leq b}(x-a)(b-x)|x-c|=\frac{(b-a)^3}{12\sqrt{3}}\,,
$
we get
\begin{equation}\label{error_interpolation_litterature}
\D|e(x)| \leq \frac{(x-a)(b-x)|x-c|}{6}\sup_{a\leq x\leq b}|f'''(x)|
\leq \D \frac{(b-a)^3}{72\sqrt{3}}\sup_{a\leq x\leq b}|f'''(x)|.
\end{equation}
This result is usually derived by considering, for $x \in [a,b]$, $x$ different from $a,b$ and $c$, the function $g(t)$ defined by
$$
g(t) = f(t)-\Pi_{[a,b]}(f)(t) - \biggl(f(x)-\Pi_{[a,b]}(f)(x)\biggr)\frac{(t-a)(t-c)(t-b)}{(x-a)(x-c)(x-b)}\,.
$$
By construction, $g(t)$ vanishes on $a, c, b$ and on the point $t=x$. Then, by applying three times Rolle's theorem, we obtain that there exists a point $\xi_x \in ]a,b[ \setminus \{c,x\}$ such that $g'''(\xi_x)=0$.\sa
Moreover, the third-order derivative of $\Pi_{[a,b]}(f)(t)$, which is a polynomial of degree 2, vanishes, and the third-order derivative of the function $t \longrightarrow (t-a)(t-c)(t-b)$ is a constant equal to $6$. Therefore, we obtain that there exists $\xi_x \in ]a,b[ \setminus \{c,x\}$ such that
$$
g'''(\xi_x)= f'''(\xi_x)- \biggl(f(x)-\Pi_{[a,b]}(f)(x)\biggr)\frac{6}{(x-a)(x-c)(x-b)}=0\,,
$$
that leads to (\ref{error_interpolation_litterature}).\\

\noindent Now, to evaluate the difference between the classical Taylor's formula and the formula we derived in (\ref{Generalized_Taylor-bis}), we consider the two constants $m_3$ and $M_3$ introduced in (\ref{m3M3}), and we reformulate estimate (\ref{error_interpolation_litterature}) by using the classical Taylor formula (\ref{m3M3}). We get the following result:
\begin{lemma}\label{error_interpolation_classique_Thm}
Let $f$ be a function $C^3([a,b])$ satisfying (\ref{m3M3}). Then the second order Taylor's theorem leads to the following interpolation error estimate:
\begin{equation}\label{error_P2_m_M}
|e(x)| \leq \frac{(b-a)^{3}}{72\sqrt{3}} (2M_{3}-m_3)\,.
\end{equation}
\end{lemma}
\begin{prooff}
We begin by writing the Lagrange $P_2-$polynomial $\Pi_{[a,b]}(f)$ given by (\ref{P2_Interpolate}) by using the classical second order Taylor's formula (\ref{Second_Order_Taylor}). \sa
For this purpose,  we substitute in (\ref{P2_Interpolate}), $f(a)$, $f(b)$ and $f(c)$ expressed in the following form:
\begin{eqnarray*}
f(a) & = & f(x) + (a-x)f'(x) + \frac{(a-x)^2}{2} f''(x) + (a-x)^2 \epsilon_{x,2}(a), \\ [0.2cm]
f(b) & = & f(x) + (b-x)f'(x) + \frac{(b-x)^2}{2} f''(x) + (b-x)^2 \epsilon_{x,2}(b), \\ [0.2cm]
f(c) & = & f(x) + (c-x)f'(x) + \frac{(c-x)^2}{2} f''(x) + (c-x)^2 \epsilon_{x,2}(c),
\end{eqnarray*}
where, by the help of (\ref{Epsilon_b_bounded}) and (\ref{m3M3}), the remainders $\epsilon_{x,2}(a)$, $\epsilon_{x,2}(b)$ and $\epsilon_{x,1}(c)$ satisfy the inequations, with $M=\max\{|m_3|,|M_3|\}$,
$$
\D|\epsilon_{x,2}(a)| \leq \frac{(x-a)}{6}M , \hs|\epsilon_{x,2}(b)| \leq \frac{(b-x)}{6}M \hs \mbox{ and } \hs |\epsilon_{x,2}(c)| \leq \frac{|c-x|}{6}M.
$$
Then, (\ref{P2_Interpolate}) gives:
\begin{eqnarray*}
\Pi_{[a,b]}(f)(x) &=&\frac{(x-c)(x-b)}{(a-c)(a-b)} \left( f(x) + (a-x)f'(x) + \frac{(a-x)^2}{2} f''(x) + (a-x)^2 \epsilon_{x,2}(a)\right) \\[0.2cm]
& + &  \frac{(x-a)(x-c)}{(b-a)(b-c)} \left( f(x) + (b-x)f'(x) + \frac{(b-x)^2}{2} f''(x) + (b-x)^2 \epsilon_{x,2}(b) \right) \\[0.2cm]
&+ &  \frac{(x-a)(x-b)}{(c-a)(c-b)} \left( f(x) + (c-x)f'(x) + \frac{(c-x)^2}{2} f''(x) + (c-x)^2 \epsilon_{x,2}(c) \right).
\end{eqnarray*}
From this expression, we can compute the coefficients of $f(x)$, $ f'(x)$ and $f''(x)$.  We get that the first one is equal to $1$, whereas the two others are equal to $0$. Consequently, we obtain for $\Pi_{[a,b]}(f)(x)$:
\begin{eqnarray}\label{Error_P2_Taylor_V0}
\D \Pi_{[a,b]}(f)(x) & = & f(x) + \frac{(x-c)(x-b)}{(a-c)(a-b)} (a-x)^2 \epsilon_{x,2}(a)
+\frac{(x-a)(x-c)}{(b-a)(b-c)} (b-x)^2 \epsilon_{x,2}(b) \nonumber \\
&& \hspace*{6.25cm}+\frac{(x-a)(x-b)}{(c-a)(c-b)} (c-x)^2  \epsilon_{x,2}(c)\,.
\end{eqnarray}
In order to determine the error $e(x)$ introduced above, we have to compute the  three last terms involved in (\ref{Error_P2_Taylor_V0}), namely
\begin{eqnarray*}
&&\hspace*{-0.8cm}\,\,\,\,\frac{(x-c)(x-b)}{(a-c)(a-b)} (a-x)^2 \epsilon_{x,2}(a)
+\frac{(x-a)(x-c)}{(b-a)(b-c)} (b-x)^2 \epsilon_{x,2}(b)
+\frac{(x-a)(x-b)}{(c-a)(c-b)} (c-x)^2 \epsilon_{x,2}(c)\,.
\end{eqnarray*}
Now recall that, from the classical Taylor formula, $\epsilon_{x,2}(a)=\ds\frac{a-x}{6}f'''(\xi(a,x))$ (the same for $\epsilon_{x,2}(b)$
and $\epsilon_{x,2}(c)$).  Using that $c$ is the midpoint of $[a,b]$, we have $\ds c-a=\frac{b-a}{2}$, and the expression above is equal to
\begin{eqnarray}\label{estimp3x}
&\hspace*{-0.6cm}  &\hspace*{-0.2cm} \frac{(a-x)(b-x)(c-x)}{3(b-a)^2}
\left[ (x-a)^{2}f'''(\xi(a,x)) + (b-x)^{2}f'''(\xi(b,x)) - 2 (c-x)^{2}f'''(\xi(c,x)) \right].
\end{eqnarray}
To estimate (\ref{estimp3x}), we will bound separately the two terms involved. For the first one,  by studying the function $f(x)=(a-x)(b-x)(c-x)$, we   easily obtain that
$$
-\frac{(b-a)}{36\sqrt{3}} \leq  \frac{(a-x)(b-x)(c-x)}{3(b-a)^2} \leq  \frac{(b-a)}{36\sqrt{3}}\,.
$$
or equivalently that
$$
\frac{|(a-x)(b-x)(c-x)|}{3(b-a)^2} \leq  \frac{(b-a)}{36\sqrt{3}}\,.
$$
\noindent For the second term of (\ref{estimp3x}), using (\ref{m3M3}), we obtain that
$$
g_{min}(x) \leq (x-a)^{2}f'''(\xi(a)) + (b-x)^{2}f'''(\xi(b)) - 2 (c-x)^{2}f'''(\xi(c)) \leq g_{max}(x)\,,
$$
where
\begin{eqnarray*}
&&g_{min}(x)= m_{3} (x-a)^{2} + m_{3} (b-x)^{2} - 2M_{3} (c-x)^{2}\,,\\
&&g_{max}(x)= M_{3} (x-a)^{2} + M_{3} (b-x)^{2} - 2m_{3} (c-x)^{2}\,.
\end{eqnarray*}
To continue, we have to determine the extremum values of the functions $g_{min}(x)$ and $g_{max}(x)$. We will consider only $g_{max}(x)$, the case of $g_{min}(x)$ being analogous. A simple computation gives that the maximum of $g_{max}$ is reached at the boundaries $x=a, x=b$ ,
with
$$
g_{max}(a)=g_{max}(b)=(b-a)^{2}\biggl(M_{3}-\ds\frac{m_3}{2}\biggr)\,.
$$
Similarly, we obtain the minimum of $g_{min}(x)$ is equal to
$$
g_{min}(a)=g_{min}(b)=(b-a)^{2}\biggl(m_{3}-\ds\frac{M_{3}}{2}\biggr)\,,
$$
so that the second term of (\ref{estimp3x}) can be bounded as follows:
$$
 (2 m_{3}-M_3)\ds\frac{(b-a)^{2}}{2} \leq (x-a)^{2}f'''(\xi(a)) + (b-x)^{2}f'''(\xi(b)) - 2 (c-x)^{2}f'''(\xi(c)) \leq (2 M_{3}-m_3) \ds\frac{(b-a)^2}{2}\,.
$$

\noindent Now, putting these results together, we finally get that
$$
|e(x)| \leq \frac{(b-a)^{3}}{72\sqrt{3}} (2M_{3}-m_3)
$$
\vspace*{-0.3cm}
\end{prooff}

\begin{remark}
This result can be compared with the more classical one recalled in (\ref{error_interpolation_litterature}). In fact the only difference comes from the term $\sup_{a\leq x\leq b}|f'''(x)|$ that is replaced here by $2M_{3}-m_3$, that takes care of the difference between the $\sup$ and the $\inf$ of the function $f'''(x)$, rather than considering their maximum.
\end{remark}

\noindent Let us now derive the corresponding result when we use the new second order Taylor-like formula (\ref{Generalized_Taylor-bis}) in the expression of the interpolation polynomial $\Pi_{[a,b]}(f)$ defined by (\ref{P2_Interpolate}). This is the purpose of the following lemma.
\begin{lemma}\label{New_error_interpolation_Thm}
Let $f\in C^3([a,b])$, then we have the following interpolation error estimate:
\begin{equation}\label{error_P1_m_M_New}
\forall x \in [a,b]: \biggl|f(x) - \Pi^*_{[a,b]}(f)(x) \biggr| \leq \frac{(b-a)^{3}}{1536 \sqrt{3}}(2M_3-m_3)\,,
\end{equation}
where $ \Pi^*_{[a,b]}(f)(x)$ is defined by
\begin{eqnarray}\label{corrected_polynomial}
&& \hspace*{-1.cm}\Pi^*_{[a,b]}(f)(x) = \Pi_{[a,b]}(f)(x) \nonumber\\
&&- \ds\frac{(x-a)(b-x)(c-x)}{(b-a)^2}\left[\ds\frac{f'(a)-2f'(c)+f'(b)}{2} + f'\biggl(\ds\frac{x+a}{2}\biggr)-2f'\biggl(\ds\frac{x+c}{2}\biggr) +f'\biggl(\ds\frac{x+b}{2}\biggr) \right ]\nonumber\\[0.1cm]
&&-\ds\frac{3 (x-a)(b-x)(c-x)}{64 (b-a)^2}\left(f''(a)(a-x) +2 f''(c)(x-c)+f''(b)(b-x)  \right)
\end{eqnarray}
\end{lemma}
\begin{prooff}
We begin to write $f(a)$, $f(b)$ and $f(c)$ by the help of (\ref{Generalized_Taylor-bis}):
$$
\D f(a) = f(x) + (a-x) \frac{f'(a) +2 f'\bigg(\ds\frac{x + a}{2}\bigg)+f'(x)}{4}
-\frac{3 (a-x)^2}{128}\left(f''(a)-f''(x)\right)
+ (a-x)^2\epsilon_{x,3}^{(2)}(a),
$$
$$
\D f(b) = f(x) + (b-x) \frac{f'(b) +2 f'\bigg(\ds\frac{x + b}{2}\bigg)+f'(x)}{4}
-\frac{3 (b-x)^2}{128}\left(f''(b)-f''(x)\right)
+ (b-x)^2\epsilon_{x,3}^{(2)}(b),
$$
$$
\D f(c) = f(x) + (c-x) \frac{f'(c) +2 f'\bigg(\ds\frac{x + c}{2}\bigg)+f'(x)}{4}
-\frac{3 (c-x)^2}{128}\left(f''(c)-f''(x)\right)
+ (c-x)^2\epsilon_{x,3}^{(2)}(c),\vspace{0.05cm}
$$
where $\epsilon_{x,3}$ satisfies (\ref{Control_Epsilon2}) with obvious changes of notations. More precisely, we have
\begin{equation}\label{3eps}
|\epsilon_{x,3}^{(2)}(a)| \leq \frac{(x-a)}{384}(2M_3-m_3)\,,\quad
|\epsilon_{x,3}^{(2)}(b)| \leq \frac{(b-x)}{384}(2M_3-m_3)\,,\quad
|\epsilon_{x,3}^{(2)}(c)| \leq \frac{|c-x|}{384}(2M_3-m_3)\,.
\end{equation}
Then, by substituting $f(a)$, $f(b)$ and $f(c)$ in the interpolation polynomial (\ref{P2_Interpolate}), we obtain
\begin{eqnarray*}
&&\hspace*{-0.8cm}\Pi_{[a,b]}(f)(x)\! =\! \\
&&\ds\frac{(x-c)(x-b)}{(a-c)(a-b)}\bigg(f(x)\! + \!(a\!-\!x) \frac{f'(a) \!+\!2 f'\biggl(\ds\frac{x + a}{2}\biggr)\!+\!f'(x)}{4}
\!-\!\frac{3 (a-x)^2}{128}\!\left(f''(a)\!-\!f''(x)\right)\!\!+\! (a-x)^2\epsilon_{x,3}^{(2)}(a)\bigg)\\
&&+ \ds\frac{(x-a)(x-c)}{(b-a)(b-c)} \bigg(f(x)\! + \!(b\!-\!x) \frac{f'(b) \!+\!2 f'\biggl(\ds\frac{x + b}{2}\biggr)\!+\!f'(x)}{4}
\!-\!\frac{3 (b-x)^2}{128}\!\left(f''(b)\!-\!f''(x)\right)\!\!+\! (b-x)^2\epsilon_{x,3}^{(2)}(b)\bigg)\\
&&+ \ds\frac{(x-a)(x-b)}{(c-a)(c-b)} \bigg(f(x)\! + \!(c\!-\!x) \frac{f'(c) \!+\!2 f'\biggl(\ds\frac{x + c}{2}\biggr)\!+\!f'(x)}{4}
\!-\!\frac{3 (c-x)^2}{128}\!\left(f''(c)\!-\!f''(x)\right)\!\!+\! (c-x)^2\epsilon_{x,3}^{(2)}(c)\bigg)
\end{eqnarray*}
From this expression, we obtain that the coefficients of $f(x)$ is equal to $1$, whereas the one before $ f'(x)$ is equal to $0$.\\

\noindent Let us compute now the terms before the other derivatives of $f(x)$. For the first order derivatives, we have the following expressions
$$
\ds\frac{(x-a)(b-x)(c-x)}{(b-a)^2}\left[\ds\frac{f'(a)-2f'(c)+f'(b)}{2} + f'\biggl(\ds\frac{x+a}{2}\biggr)-2f'\biggl(\ds\frac{x+c}{2}\biggr) +f'\biggl(\ds\frac{x+b}{2}\biggr) \right ]\,,
$$
whereas the terms in $f''(x)$ can be expressed as
$$
\ds\frac{3 (x-a)(b-x)(c-x)}{64 (b-a)^2}\left(f''(a)(a-x) +2 f''(c)(x-c)+f''(b)(b-x)  \right)\,.
$$
Finally, we obtain for $\Pi_{[a,b]}(f)(x)$:
\begin{eqnarray}\label{newPoly}
&&\hspace*{-0.7cm}\Pi_{[a,b]}(f)(x)\! = f(x)\!
+\ds\frac{(x-a)(b-x)(c-x)}{(b-a)^2}\bigg[\ds\frac{f'(a)-2f'(c)+f'(b)}{2} + f'\biggl(\ds\frac{x+a}{2}\biggr)-2f'\biggl(\ds\frac{x+c}{2}\biggr) +f'\biggl(\ds\frac{x+b}{2}\biggr) \bigg]\nonumber\\
&&\hspace*{0.5cm}+\ds\frac{3 (x-a)(b-x)(c-x)}{64 (b-a)^2}\bigg(f''(a)(a-x) +2 f''(c)(x-c)+f''(b)(b-x)  \bigg)\nonumber\\
&&\hspace*{0.5cm}+\ds\frac{(x-c)(x-b)}{(a-c)(a-b)}  (a-x)^2\epsilon_{x,3}^{(2)}(a)
+\ds\frac{(x-a)(x-c)}{(b-a)(b-c)}  (b-x)^2\epsilon_{x,3}^{(2)}(b)
+\ds\frac{(x-a)(x-b)}{(c-a)(c-b)}(c-x)^2\epsilon_{x,3}^{(2)}(c).
\end{eqnarray}
Now, let us consider the new interpolation polynomial $\Pi^*_{[a,b]}(f)$ introduced in (\ref{corrected_polynomial}). With this, the interpolation polynomial $\Pi_{[a,b]}(f)$ of (\ref{newPoly}) can be expressed as
$$
\Pi^*_{[a,b]}(f) = f(x)+E(x,a,b)\,,
$$
where the function error $E(x,a,b)$ is defined by
 $$
E(x,a,b)= \ds\frac{(x-c)(x-b)}{(a-c)(a-b)}  (a-x)^2\epsilon_{x,3}^{(2)}(a)
+\ds\frac{(x-a)(x-c)}{(b-a)(b-c)}  (b-x)^2\epsilon_{x,3}^{(2)}(b)
+\ds\frac{(x-a)(x-b)}{(c-a)(c-b)}(c-x)^2\epsilon_{x,3}^{(2)}(c)\,.
$$
Using (\ref{3eps}), we obtain the following bound
\begin{equation}\label{useone}
|E(x,a,b)| \leq 2 \frac{|(a-x)(b-x)(c-x)|}{(b-a)^2} \frac{2M_3-m_3}{384} \left( (x-a)^{2} + (b-x)^{2} +2 (c-x)^{2} \right)\,.
\end{equation}
We already saw before that
$$
\ds\frac{|(a-x)(b-x)(c-x)|}{(b-a)^2} \leq \frac{b-a}{12 \sqrt{3}}\,,
$$
so we only have to bound the second term of (\ref{useone}), that is $(x-a)^{2} + (b-x)^{2} +2 (c-x)^{2}$.
\sa
Using a method similar to the one used for (\ref{estimp3x}), we obtain that
$$
\max_{a \leq x \leq b} \left( (x-a)^{2} + (b-x)^{2} +2 (c-x)^{2} \right)= \frac{3}{2}(b-a)^{2}\,.
$$
Putting all together, we obtain that
$$
|E(x,a,b)| \leq \frac{2M_3-m_3}{1536 \sqrt{3}}(b-a)^{3}\,,
$$
which completes the proof of this lemma.
\end{prooff}

\begin{remark}
Let us compare the interpolation errors of lemma \ref{New_error_interpolation_Thm} with the classical one of lemma  \ref{error_interpolation_classique_Thm}. First, note that $ \Pi^*_{[a,b]}(f)$ is a polynomial of degree less than or equal to 3. So, interpolation errors (\ref{error_P2_m_M}) and (\ref{error_P1_m_M_New}) can not be compared anymore, since $ \Pi_{[a,b]}(f)$ is a polynomial of degree less than or equal to 2.
\sa
However, if we consider the Lagrange polynomial $ \Pi^{(3)}_{[a,b]}(f)$ of degree less than or equal to 3, we can compare the corresponding interpolation error with  result (\ref{error_P1_m_M_New}), assuming that the function $f \in C^{4}([a,b])$. Indeed, the standard estimate in literature requires the function $f$ to belong to $ C^{4}([a,b])$.
\sa
Hence, denoting by $M_4=\ds\sup_{a\leq x\leq b}f^{(4)}(x)$ and by $m_4=\ds\inf_{a\leq x\leq b}f^{(4)}(x)$, the Lagrange interpolating polynomial error bound can be written (see for instance \cite{Atki88}, \cite{SuMa03})
\begin{equation}\label{error_Pi3}
  \biggl|f(x) - \Pi^{(3)}_{[a,b]}(f)(x) \biggr| \leq \frac{(b-a)^{4}}{1296}{\cal M}_4, \mbox{ with }{\cal M}_4:=\max(|M_4|,|m_4|)\,.
\end{equation}
 Now, to compare this estimate with (\ref{error_P1_m_M_New}), we assume that variations of the third derivative $f^{(3)}$ are not ``extreme'', so that $2M_{3}-m_3$ can approximatively be replaced  by ${\cal M}_{3}:=\max(|m_3|, |M_3|)$. In addition, for a function $f$ smooth enough (at least $C^{4}([a,b])$), we can (roughly speaking) also assume that $|f^{(3)}|$ behaves like $(b-a)|f^{(4)}|$. Hence, the error estimate  (\ref{error_P1_m_M_New}) can be expressed as
 $$
  \biggl|f(x) - \Pi^*_{[a,b]}(f)(x) \biggr| \leq \frac{(b-a)^{4}}{1536 \sqrt{3}}{\cal M}_{4}\,.
 $$
This gives a bound which is more than 2 times smaller than the corresponding one in (\ref{error_Pi3}).
\end{remark}
\begin{remark}
 Let us give an elementary numerical example. Consider the interval $[a,b]=[0,1]$, and the function $f(x)=\ln(1+x)$. Formula (\ref{Generalized_Taylor-bis}) gives
$$
\ln(2)=\ln(1)+\ds\frac{\frac{1}{2}+2\frac{2}{3}+1}{4}-\frac{3}{128}(-\frac{1}{4}+1)+\epsilon_{a,3}^{(2)} =\ds\frac{1061}{1536} +\epsilon_{a,3}^{(2)}\simeq 0.6907+\epsilon_{a,3}^{(2)}
$$
and
$$
|\epsilon_{a,3}^{(2)}| \leq \ds\frac{15}{1536} \simeq 1/100\,.
$$
With the same data,  classical Taylor's formula (\ref{Second_Order_Taylor}) gives
$$
\ln(2)=\ln(1)+1+\frac{1}{2}(-1)+\epsilon_{a,2} = \frac{1}{2}+\epsilon_{a,2}
$$
and
$$
|\epsilon_{a,2}| \leq \ds\frac{1}{3}\,.
$$
Hence, the improved formula leads to a much more accurate approximation of $\ln(2)$.
\end{remark}
\subsection{The quadrature error}\label{Quadrature_Error}
\noindent We consider now, for any integrable function $f$ defined on $[a,b]$, Simpson's quadrature rule \cite{Crouzeix_Mignot} whose formula is given by
\begin{equation}\label{Simpson}
\D  \int_{a}^{b} f(x)dx \simeq \frac{b-a}{6}\left(f(a)+4\,f\left(\frac{a+b}{2}\right)+f(b)\right).
\end{equation}
The reason we consider (\ref{Simpson}) is that this quadrature formula corresponds to approximate the function $f$ by its Lagrange polynomial interpolation $\Pi_{[a,b]}(f)$, of degree less than or equal to two, which is given by (\ref{P2_Interpolate}).
\sa
\noindent Thus, in the classical literature of numerical integration (see for example \cite{DrACxx}, \cite{Crouzeix_Mignot} and \cite{Cerone}), we can find the standard \emph{Simpson inequality}
\begin{equation}\label{Erreur_Simpson_Standard}
\D\left|\int_{a}^{b}f(x)\,dx -  \frac{b-a}{6}\left(f(a)+4\,f\left(\frac{a+b}{2}\right)+f(b)\right)\right| \leq \frac{(b-a)^5}{2880}\sup_{a\leq x\leq b}|f^{(4)}(x)|,
\end{equation}
for any function four times differentiable $f$ on $[a,b]$, whose fourth derivative is accordingly bounded on $[a,b]$.
\sa
\noindent Now, if the function $f$ is not four times differentiable, or if the fourth derivative $f^{(4)}$ is not bounded on $[a,b]$, we cannot apply the formula above. Therefore, if we consider a function $f$ that is only $\mathcal{C}^{3}$ on $[a,b]$, we have the following estimate \cite{BuFa11}, \cite{ChSu02}
\begin{equation}\label{Erreur_Simpson_C1}
\D\left|\int_{a}^{b}f(x)\,dx -  \frac{b-a}{6}\left(f(a)+4\,f\left(\frac{a+b}{2}\right)+f(b)\right)\right| \leq \frac{(b-a)^4}{192}\sup_{a\leq x\leq b}|f^{'''}(x)|,
\end{equation}
Let us prove now a result that gives estimate (\ref{Erreur_Simpson_C1}) in an alternative display, based on the classical Taylor formula.
\begin{lemma}\label{lemma33}
Let $f \in C^{3}([a,b])$ which satisfies (\ref{m3M3}). Then, we have the following estimate:
$$
\D\left|\int_{a}^{b}f(x)\,dx -  \frac{b-a}{6}\left(f(a)+4\,f\left(\frac{a+b}{2}\right)+f(b)\right)\right| \leq  \frac{5 (b-a)^4}{1152}(M_{3}-m_{3})\,.
$$
\end{lemma}
\begin{prooff}
To derive this estimate, let us begin with the classical second order Taylor's formula (\ref{Second_Order_Taylor}), from which we have derived above the expression (\ref{Error_P2_Taylor_V0}) of polynomial $\Pi_{[a,b]}(f)$. Then, by integrating the difference $f(x) - \Pi_{[a,b]}(f) (x)$ between $a$ and $b$, we get
\begin{eqnarray}
\label{Erreur_Simpson_C2_V01}
\hspace*{-0.9cm}\int_{a}^{b}\!\! (f(x) - \Pi_{[a,b]}(f) (x)) \,dx \!\!\!&= &\!\!\!\! \int_{a}^{b} \left( \frac{(x-c)(b-x)}{(c-a)(b-a)} \frac{(a-x)^2}{2}\epsilon_{x,2}(a) \right.\nonumber \\[0.2cm]
&+&\!\!\!\! \!\! \left. \frac{(x-a)(c-x)}{(b-a)(b-c)}\frac{(b-x)^2}{2}\epsilon_{x,2}(b)\! +\!\frac{(x-a)(x-b)}{(c-a)(b-c)} \frac{(c-x)^2}{2}\epsilon_{x,2}(c) \!  \right) \!dx.
\end{eqnarray}

\noindent However, it is well known \cite{Atki88} that the $P_2-$Lagrange interpolation polynomial $\Pi_{[a,b]}(f)$ given by (\ref{P2_Interpolate}) also fulfills:
\begin{equation}\label{integrale_P2_V0}
\int_{a}^{b}\Pi_{[a,b]}(f)(x)\,dx =  \frac{b-a}{6}\left(f(a)+4 f(c)+ f(b) \right).
\end{equation}
Now, let us introduce the usual error in the quadrature rule $E(f)$ defined by
$$
E(f) \equiv \int_{a}^{b} f(x)dx - \frac{b-a}{6}\left(f(a)+4 f(c)+ f(b) \right)\,.
$$
Using the expressions $\epsilon_{x,2}(a), \epsilon_{x,2}(b)$ and $\epsilon_{x,2}(c)$ (see proof of lemma \ref{error_interpolation_classique_Thm}), equations (\ref{Erreur_Simpson_C2_V01}) and (\ref{integrale_P2_V0}) give
\begin{eqnarray*}
E(f) &= & \int_{a}^{b} \frac{(x-a)(b-x)(c-x)}{3(b-a)^2} \left( (x-a)^2 f'''(\xi(a)) + (b-x)^2 f'''(\xi(b))  -2 \, (c-x)^2 f'''(\xi(c))  \,  \right) \!dx
\end{eqnarray*}
that we split, for convenience, in three integrals $I(a), I(b)$ and $I(c)$, so that
$$
E(f)=I(a) + I(b)- I(c)
\vspace{-0.3cm}
$$
with
\begin{eqnarray*}
&&I(a) =  \frac{1}{3(b-a)^2} \int_{a}^{b} (x-a)^3(b-x)(c-x) f'''(\xi(a)) dx\,, \\
&&I(b) =  \frac{1}{3(b-a)^2} \int_{a}^{b} (x-a)(b-x)^3(c-x) f'''(\xi(b)) dx \,,\\
&&I(c) =  \frac{2}{3(b-a)^2} \int_{a}^{b} (x-a)(b-x)(c-x)^3 f'''(\xi(c)) dx\,.
\end{eqnarray*}
To obtain estimates of $E(f)$, we will consider separately each of these integrals. We detailed here the computations for $I(a)$, the others terms can be treated similarly.
\sa
\noindent Noting that the term $(x-a)^3(b-x)(c-x)$ is positive for $a \leq x \leq c$ and negative for $c \leq x \leq b$, we split the integral and we get, using the mean value theorem, that there exists a constant $a < C_{1,a} < c$, (respectively $c < C_{2,a} < b$), such that
\begin{eqnarray*}
&& \int_{a}^{c} (x-a)^3(b-x)(c-x) f'''(\xi(a)) dx =  f'''(C_{1,a}) \int_{a}^{c} (x-a)^3(b-x)(c-x) dx \,,\\
&& \int_{c}^{b} (x-a)^3(b-x)(c-x) f'''(\xi(a)) dx =  f'''(C_{2,a}) \int_{c}^{b} (x-a)^3(b-x)(c-x) dx\,.
\end{eqnarray*}
It remains now to compute the integral above. Owing to the relation $c=(a+b)/2$, a straightforward computation gives us
$$
\int_{a}^{c} (x-a)^3(b-x)(c-x) dx = \frac{(b-a)^6}{960}, \quad \mbox{ and }  \int_{c}^{b} (x-a)^3(b-x)(c-x) dx = -\frac{3(b-a)^6}{320}\,.
$$
Using now inequalities (\ref{m3M3}), we readily get
\begin{eqnarray*}
&& (m_{3}-9 M_{3}) \frac{(b-a)^4}{2880} \leq I(a) = \frac{(b-a)^4}{2880} \biggl(f'''(C_{1,a})- 9 f'''(C_{2,a})\biggr) \leq (M_{3}-9m_{3})\frac{(b-a)^4}{2880}\,.
\end{eqnarray*}

\noindent The same computations for $I(b)$  gives
$$
\int_{a}^{c} (x-a)(b-x)^3(c-x) dx =  \frac{3(b-a)^6}{320} , \quad \mbox{ and }  \int_{c}^{b} (x-a)(b-x)^3(c-x) dx = -\frac{(b-a)^6}{960}\,,
$$
so that
\begin{eqnarray*}
&& (9m_{3}- M_{3}) \frac{(b-a)^4}{2880} \leq I(b) = \frac{(b-a)^4}{2880} \biggl(9 f'''(C_{1,b})- f'''(C_{2,b})\biggr) \leq  (9 M_{3}-m_{3})\frac{(b-a)^4}{2880}\,.
\end{eqnarray*}

\noindent Similarly for $I(c)$, where the polynomial in the integral is odd with respect to $x=c$, we get
$$
\int_{a}^{c} (x-a)(b-x)(c-x)^3 dx =  \frac{(b-a)^6}{768} , \quad \mbox{ and }  \int_{c}^{b} (x-a)(b-x)(c-x)^3 dx = -\frac{(b-a)^6}{768}\,,
$$
so that
\begin{eqnarray*}
&& (m_{3}- M_{3}) \frac{(b-a)^4}{1152} \leq I(c) = \frac{(b-a)^4}{1152} \biggl(f'''(C_{1,c})- f'''(C_{2,c})\biggr) \leq (M_{3}-m_{3})\frac{(b-a)^4}{1152}\,.
\end{eqnarray*}
Putting all together, we finally proved that
\begin{eqnarray*}
&& 5 (m_{3}- M_{3}) \frac{(b-a)^4}{1152} \leq E(f) = I(a)+I(b)-I(c) \leq 5 (M_{3}-m_{3})\frac{(b-a)^4}{1152}\,.
\end{eqnarray*}
\end{prooff}
\noindent Let us consider now the corresponding quadrature formula denoted $I(f)$, based on the generalized Taylor-like expansion (\ref{Generalized_Taylor-bis}) and defined by:
\begin{eqnarray}\label{formula-improved}
\hspace*{-0.5cm} I(f)  & = &  \frac{b-a}{6}\bigg[f(a)+4\,f\bigg(\frac{a+b}{2}\bigg)+f(b)\bigg] \nonumber \\
&& \nonumber  \\
& + &\hspace*{-0.4cm} \frac{(b-a)^3}{240} \left[f''\biggl(\ds \frac{a}{2} \biggr)  -2 f''\biggl(\ds \frac{c}{2} \biggr) + f''\biggl(\ds \frac{b}{2} \biggr) \right]
 - \frac{(b-a)^3}{2560}\big(f''(a)+2 f''(c)+f''(b)\big)
\end{eqnarray}
The corresponding quadrature error estimate is then given by the following lemma:
\begin{lemma}\label{New_Quadrature_Error_Thm}
Let $f \in C^3([a,b])$ which satisfies (\ref{m3M3}), with a third derivative $f'''$ L-Lipschitz. Then, we have the following estimate:
\begin{equation}\label{new formulaQuadra}
\D\left| \int_{a}^{b}f(x)\,dx -  I(f) \right| \,\leq \,\frac{L \,(b-a)^3(a^2+ab+b^2)}{512}+\frac{5(b-a)^4}{36864}(2M_3-m_3)\,.
\end{equation}
Moreover, if $ab < 0$, we have:
\begin{equation}\label{new formulaQuadraBIS}
\D\left| \int_{a}^{b}f(x)\,dx -  I(f) \right| \,\leq \,\frac{L \,(b-a)^5}{512}+\frac{5(b-a)^4}{36864}(2M_3-m_3)\,.
\end{equation}
\end{lemma}
\begin{prooff}
We consider the expression (\ref{newPoly}) of the polynomial function $\Pi_{[a,b]}(f)(x)$. By integrating between $a$ and $b$ the difference $f(x)- \Pi_{[a,b]}(f)(x)$, we obtain:
\begin{eqnarray}
&&\hspace*{-0.5cm}\D\int_{a}^{b}\hspace{-0.3cm}\left(f(x)- \Pi_{[a,b]}(f)(x)\right)dx =\nonumber\\[0.2cm]
&&\hspace*{-0.5cm}- \D\int_{a}^{b} \ds\frac{(x-a)(b-x)(c-x)}{(b-a)^2}\left[\ds\frac{f'(a)-2f'(c)+f'(b)}{2} + f'\bigg(\ds\frac{x+a}{2}\bigg)-2f'\bigg(\ds\frac{x+c}{2}\bigg) +f'\bigg(\ds\frac{x+b}{2}\bigg) \right]\label{clear1}\\[0.2cm]
&&\hspace*{-0.5cm}-\D\int_{a}^{b}\ds\frac{3 (x-a)(b-x)(c-x)}{64 (b-a)^2}\left(f''(a)(a-x) +2 f''(c)(x-c)+f''(b)(b-x)  \right) dx\label{clear2}\\[0.2cm]
&&\hspace*{-0.5cm}-\D\int_{a}^{b} \ds\frac{(x-c)(x-b)}{(a-c)(a-b)}  (a-x)^2\epsilon_{x,3}^{(2)}(a)
+\ds\frac{(x-a)(x-c)}{(b-a)(b-c)}  (b-x)^2\epsilon_{x,3}^{(2)}(b)
+\ds\frac{(x-a)(x-b)}{(c-a)(c-b)}(c-x)^2\epsilon_{x,3}^{(2)}(c) dx\nonumber \\[0.2cm] \label{clear3}
\end{eqnarray}
In the following of the proof, we will consider one by one, each line (\ref{clear1}), (\ref{clear2}) and (\ref{clear3}) of the above formula:\\
\noindent Given that
$$
\D\int_{a}^{b}(x-a)(b-x)(c-x) dx =0\,,
$$
the first line (\ref{clear1}) can be written as
\begin{equation}\label{line1}
- \frac{1}{(b-a)^2}\D\int_{a}^{b} (x-a)(b-x)(c-x) F'_{a,b,c}(x) dx\,,
\end{equation}
where we denote $F'_{a,b,c}(x)= f'\bigg(\ds\frac{x+a}{2}\bigg)-2f'\bigg(\ds\frac{x+c}{2}\bigg) +f'\bigg(\ds\frac{x+b}{2}\bigg)$.
\sa
We begin to write $\ds f'\biggl(\ds \frac{x+a}{2}\biggr)$, $\ds f'\biggl(\ds \frac{x+c}{2}\biggr)$ and $\ds f'\biggl(\ds \frac{x+b}{2}\biggr)$,  by the help of the first order Taylor's formula. We get
$$
\ds f'\bigg(\ds\frac{x+a}{2}\bigg)= \ds f'\biggl(\ds \frac{a}{2}\biggr)+ \frac{x}{2} \ds f''\biggl(\ds \frac{a}{2}\biggr) + \frac{x^2}{8} \ds f'''\bigl(\xi_1(x)\bigr)\,, \quad \ds \frac{a}{2}<  \xi_1(x) < \ds\frac{x+a}{2}\,,
$$
$$
\ds f'\bigg(\ds\frac{x+c}{2}\bigg)= \ds f'\biggl(\ds \frac{c}{2}\biggr)+ \frac{x}{2} \ds f''\biggl(\ds \frac{c}{2}\biggr) + \frac{x^2}{8} \ds f'''\bigl(\xi_2(x)\bigr)\,, \quad \ds \frac{c}{2}<  \xi_1(x) < \ds\frac{x+c}{2}\,,
$$
$$
\ds f'\bigg(\ds\frac{x+b}{2}\bigg)= \ds f'\biggl(\ds \frac{b}{2}\biggr)+ \frac{x}{2} \ds f''\biggl(\ds \frac{b}{2}\biggr) + \frac{x^2}{8} \ds f'''\bigl(\xi_3(x)\bigr)\,,  \quad \ds \frac{b}{2}<  \xi_1(x) < \ds\frac{x+b}{2}\,.
$$
Substituting these expressions in (\ref{line1}), we obtain that the term with  the first derivatives is equal to zero. The first non vanishing term in (\ref{line1}) (with the second derivatives) is equal to
\begin{eqnarray}\label{contribnew1}
&-& \frac{1}{(b-a)^2}\left(f''\biggl(\ds \frac{a}{2} \biggr)  -2 f''\biggl(\ds \frac{c}{2} \biggr) + f''\biggl(\ds \frac{b}{2} \biggr) \right) \D\int_{a}^{b} \frac{x}{2}(x-a)(b-x)(c-x)  dx \nonumber \\
&& \nonumber\\
&=&- \frac{(b-a)^3}{240} \left(f''\biggl(\ds \frac{a}{2} \biggr)  -2 f''\biggl(\ds \frac{c}{2} \biggr) + f''\biggl(\ds \frac{b}{2} \biggr) \right) \,.
\end{eqnarray}
Now, the last term of (\ref{line1}), that contributes to the error bound gives
\begin{equation}\label{LastTerm}
- \frac{1}{8(b-a)^2} \D\int_{a}^{b} \left[f'''\bigl(\xi_1(x)\bigr)  -2 f'''\bigl(\xi_2(x)\bigr) + f'''\bigl(\xi_3(x)\bigr) \right] x^2(x-a)(b-x)(c-x)  dx\,.
\end{equation}
Let us consider the term with the third derivatives of (\ref{LastTerm}). Using that $f'''$ is L-Lipschitz, we have
\begin{eqnarray*}
|f'''\bigl(\xi_1(x)\bigr) \! - \! 2 f'''\bigl(\xi_2(x)\bigr) \!+\! f'''\bigl(\xi_3(x)\bigr)| & \leq  &
|f'''\bigl(\xi_1(x)\bigr) \!- \! f'''\bigl(\xi_2(x)\bigr)| +| f'''\bigl(\xi_3(x)\bigr) \!-\!  f'''\bigl(\xi_2(x)\bigr)| \\
&&\\
& \leq & L\,|\xi_1(x)-\xi_2(x)| +  L\,|\xi_3(x)-\xi_2(x)| \\
&&\\
& \leq & \ds\frac{3}{4}L\,(b-a)+\ds\frac{3}{4}L\,(b-a) = \ds\frac{3}{2}L\,(b-a)\,,
\end{eqnarray*}
using that $a<\xi_1(x)<c$, $\ds\frac{a+c}{2}<\xi_2(x)<\ds\frac{b+c}{2}$ and $c<\xi_3(x)<b$. In these conditions, (\ref{LastTerm}) can be bounded in absolute value by
$$
 \frac{3 \, L}{16(b-a)} \D\int_{a}^{b}x^2(x-a)(b-x)|c-x| dx\,.
$$
Given that
$$
\D\int_{a}^{b}x^2(x-a)(b-x)|c-x| dx= \frac{(b-a)^4}{96}(a^2+ab+b^2)\,,
$$
(\ref{LastTerm}) can be bounded by
$$
\frac{3 \, L (b-a)^3}{1536}(a^2+ab+b^2)=\frac{L\, (b-a)^3}{512}(a^2+ab+b^2)\,.
$$
Assuming (for instance) that $ab <0$, as for example for $[a,b]=[-1,1]$, we readily get that $a^2+ab+b^2 \leq (b-a)^{2}$ and
\begin{equation}\label{boundLip}
\frac{L\, (b-a)^3}{512}(a^2+ab+b^2) \leq \frac{L\, (b-a)^5}{512}\,.
\end{equation}
\sa
\noindent Consider now the second line (\ref{clear2}) that can be decomposed into three terms. Let us consider for example the first one, namely
$$
-\frac{3f''(a)}{64 (b-a)^2}\D\int_{a}^{b} (x-a)(b-x)(c-x)(a-x) dx\,.
$$
We can explicitly compute this integral that yields
$$
\D\int_{a}^{b} (x-a)(b-x)(c-x)(a-x) dx= -\frac{(b-a)^5}{120}\,.
$$
Consequently, the first term of (\ref{clear2}) is equal to
$$
-\frac{3f''(a)}{64 (b-a)^2}\D\int_{a}^{b} (x-a)(b-x)(c-x)(a-x) dx =\frac{(b-a)^3 f''(a)}{2560} \,.
$$
Similarly, for the second term of (\ref{clear2}), we get
$$
-\frac{6f''(c)}{64 (b-a)^2}\D\int_{a}^{b} (x-a)(b-x)(c-x)(x-c) dx =
 \frac{(b-a)^3 f''(c)}{1280}\,,
$$
and, for the last one of (\ref{clear2})
$$
-\frac{3f''(b)}{64 (b-a)^2}\D\int_{a}^{b} (x-a)(b-x)(c-x)(b-x) dx\ =
\frac{(b-a)^3 f''(b)}{2560}\,.
$$
Summing up, the second line (\ref{clear2}) gives the following contribution
\begin{equation}\label{contribnew2}
\frac{(b-a)^3}{2560}\big(f''(a)+2 f''(c)+f''(b)\big) \,.
\end{equation}
\noindent Let us consider now the last line (\ref{clear3}) that can be expressed as, $c$ being the midpoint of $[a,b]$:
$$
\frac{2}{(b-a)^2}\D\int_{a}^{b}\!\! \!-(x-a)^2(b-x)(c-x) \epsilon_{x,3}^{(2)}(a)
\,+\,\ds(x-a) (b-x)^2(c-x) \epsilon_{x,3}^{(2)}(b)
\,-\,2(x-a)(b-x)(c-x)^2\epsilon_{x,3}^{(2)}(c) \,dx,
$$
the absolute value of which being bounded by the sum of three terms $J(a)+J(b)+2J(c)$, where we define
\begin{eqnarray*}
&&J(a)=\frac{2}{(b-a)^2}\D\int_{a}^{b} (x-a)^2(b-x)\big|c-x\big| \big|\epsilon_{x,3}^{(2)}(a)\big| dx\,,\\[0.1cm]
&&J(b) = \frac{2}{(b-a)^2}\D\int_{a}^{b} \ds(x-a) (b-x)^2\big|c-x\big| \big|\epsilon_{x,3}^{(2)}(b)\big| dx\,, \\[0.1cm]
&&J(c)=\frac{2}{(b-a)^2}\D\int_{a}^{b}  (x-a)(b-x)(c-x)^2\big|\epsilon_{x,3}^{(2)}(c)\big| dx\,.
\end{eqnarray*}
\noindent Now, we will bound separately each of these 3 terms, using again the estimates (\ref{3eps}) to bound $\epsilon_{x,3}^{(2)}(a), \epsilon_{x,3}^{(2)}(c)$ and $\epsilon_{x,3}^{(2)}(b)$. We  have:
\begin{eqnarray*}
&&|J(a)| \leq \frac{2}{(b-a)^2}\frac{2M_3-m_3}{384}\D\int_{a}^{b} (x-a)^3(b-x)|c-x| dx\,,\\[0.1cm]
&& |J(b)| \leq \frac{2}{(b-a)^2}\frac{2M_3-m_3}{384} \D\int_{a}^{b} (x-a)(b-x)^3|c-x| dx\,,\\[0.1cm]
&& |J(c)| \leq \frac{2}{(b-a)^2}\frac{2M_3-m_3}{384} \D\int_{a}^{b} (x-a)(b-x)|c-x|^3 dx\,.
\end{eqnarray*}
Using similar computations as above (see proof of lemma \ref{lemma33}), straightforward computations give
$$
\int_{a}^{b} (x-a)^3(b-x)|c-x| dx =  \D\int_{a}^{b} (x-a)(b-x)^3|c-x| dx= \frac{(b-a)^6}{96},
$$
and
$$
\int_{a}^{b} (x-a)(b-x)|c-x|^{3} dx = \frac{(b-a)^6}{384}\,.
$$
Putting all together, we finally obtain that the last line (\ref{clear3}) is bounded by the term
\begin{equation}\label{Tterm3}
\frac{5}{32}\frac{(b-a)^4}{1152}(2M_3-m_3)\,.
\end{equation}
Finally, combining  estimates (\ref{boundLip}) and (\ref{Tterm3}), and adding the contributions (\ref{contribnew1}) and (\ref{contribnew2}) to the quadrature approximation of the integral, we get the error estimate (\ref{new formulaQuadraBIS}).
\end{prooff}
\noindent Let us conclude this section by several remarks.
\begin{enumerate}
\item To obtain estimate (\ref{new formulaQuadraBIS}), we used that $a^2+ab+b^2 \leq (b-a)^2$. In fact, in several cases, we can derive a better estimate of the form $a^2+ab+b^2 \leq \frac{(b-a)^2}{m}, m>1$. For instance, for $[a,b]=[-1,1]$, we readily get $m=4$, leading to an improved error bound.
\item We are now interested to compare the relative numerical weights between the two terms involved in estimate (\ref{new formulaQuadraBIS}), namely
$$
\frac{L \,(b-a)^5}{512} \mbox{ and }\frac{5(b-a)^4}{36864}(2M_3-m_3).
$$
As this comparison depends on the constants $m_3$ , $M_3$ and $L$, it is not reachable in a general case. Hence, we will consider a numerical example, where we will be able to avoid the numerical evaluation of $L$.
\sa
Namely, for a given interval $[a,b]$, $(a < b$ and $ab<0)$, we introduce the function
$$
f(x)=(x-a)^p, 3 \leq p < 4, \mbox{ where } f'''(x)=p(p-1)(p-2)(x-a)^{p-3}.
$$
For these values of $p$, the function $f$ belongs to $C^3([a,b])$ but is not four times derivable at point $a$.
\sa
Then, to estimate the quadrature error (\ref{new formulaQuadraBIS}), we have to evaluate the term (\ref{LastTerm}), that is the absolute value of the quantity $X$ defined by
$$
X:=- \frac{1}{8(b-a)^2} \D\int_{a}^{b} \left[f'''\bigl(\xi_1(x)\bigr)  -2 f'''\bigl(\xi_2(x)\bigr) + f'''\bigl(\xi_3(x)\bigr) \right] x^2(x-a)(b-x)(c-x)  dx\,.
$$
As in the proof of lemma \ref{New_Quadrature_Error_Thm}, we have to bound the terms
$|f'''\bigl(\xi_1(x)\bigr) \!- \! f'''\bigl(\xi_2(x)\bigr)|$ and $| f'''\bigl(\xi_3(x)\bigr) \!-\!  f'''\bigl(\xi_2(x)\bigr)|$.
\sa
 Let us deal with the first one, the second being similar.  We have, assuming for instance $\xi_1(x)<\xi_2(x)$,
$$
f'''\bigl(\xi_1(x)\bigr) \!- \! f'''\bigl(\xi_2(x)\bigr) = p(p-1)(p-2)\bigl( (\xi_1(x)-a)^{p-3} - (\xi_2(x)-a)^{p-3} \bigr)\,.
$$
Now, recall that  $a<\xi_1(x)<c$ whereas $\ds\frac{a+c}{2}<\xi_2(x)<\ds\frac{b+c}{2}$, we obtain that
$$
-\ds\biggl(\frac{3}{4}\biggr)^{p-3}(b-a)^{p-3}\leq (\xi_1(x)-a)^{p-3} - (\xi_2(x)-a)^{p-3} \leq \ds\biggl[\biggl(\frac{1}{2}\biggr)^{p-3}- \biggl(\frac{1}{4}\biggr)^{p-3}\biggr](b-a)^{p-3}\,,
$$
so that
$$
|f'''\bigl(\xi_1(x)\bigr) \!- \! f'''\bigl(\xi_2(x)\bigr)| \leq \biggl(\frac{3}{4}\biggr)^{p-3}(b-a)^{p-3}\,.
$$
Using the same estimate for $| f'''\bigl(\xi_3(x)\bigr) \!-\!  f'''\bigl(\xi_2(x)\bigr)|$, and continuing the computation exactly as in the proof of Lemma \ref{New_Quadrature_Error_Thm}, we obtain that $|X|$ can be bounded by
$$
|X| \leq \biggl(\frac{3}{4}\biggr)^{p-3}\,\frac{p(p-1)(p-2)}{384}(b-a)^{p+1}\,,
$$
that avoids the evaluation of the Lipschitz constant $L$.
\sa
Taking into account that $m_3=0$ and $M_3=p(p-1)(p-2)(b-a)^{p-3}$, we also have
$$
\frac{5(b-a)^4}{36864}(2M_3-m_3)\, \leq \,\frac{p(p-1)(p-2)}{3686}(b-a)^{p+1} \,.
$$
Now, given that $3 \leq p <4$, we obtain
$$
\frac{1}{3686}\,<< \frac{1}{512}\, \leq \,\frac{1}{384}\biggl(\frac{3}{4}\biggr)^{p-3}\!\!\leq \,\frac{1}{384}\,,
$$
which shows that the term $\ds\frac{5(b-a)^4}{36864}(2M_3-m_3)$ is negligeable compared to the term $\ds\biggl(\frac{3}{4}\biggr)^{p-3}\,\frac{p(p-1)(p-2)}{384}(b-a)^{p+1}$.\vspace{0.1cm}
\item Interestingly, error estimate (\ref{new formulaQuadraBIS}) can also be compared to the optimized result obtained in \cite{ChSu02}, \cite{DrACxx} for a third-order differentiable function on $[a,b]$. There, it is proved that, for another numerical integration rule, we have
\begin{equation}
\label{Cheng-Sun}
\D\left|\int_{a}^{b}f(x)\,dx -  \frac{b-a}{2}\left(f(a)+f(b)\right)+\frac{(b-a)^2}{12}(f'(b)-f'(a))\right| \leq \frac{(b-a)^4}{384}(M_3-m_3)\,.
\end{equation}
To be able to compare (\ref{Cheng-Sun}) to our result, we assume that, roughly speaking, Lipschitz constant $L$ can be approximated by $L\simeq\ds\frac{M_3-m_3}{b-a}$. Hence, error bound (\ref{new formulaQuadraBIS}) can be  written
$$
\frac{(M_3-m_3) \, (b-a)^4}{512}+\frac{5(b-a)^4}{36864}(2M_3-m_3)\,.
$$
Consequently, the first term appears to be much greater than second one. Hence, in that case, the error bound can be approximated by
\begin{equation}\label{Last formulaQuadra}
\frac{ (b-a)^4}{512}(M_3-m_3)\,,
\end{equation}
and (\ref{Last formulaQuadra}) is $1.33$ times smaller than (\ref{Cheng-Sun}). However, the ``price to pay'' in the quadrature formula  (\ref{formula-improved}) is the computation of $f''$ at the points $a, \ds\frac{a+b}{2}$ and $b$, compared with $f'$ at the points $a$ and $b$
in formula (\ref{Cheng-Sun}).
\end{enumerate}
\section{Conclusion}\label{D}
\noindent In this paper we proposed a new second-order Taylor-like theorem to obtain some minimized remainders. For a function $f$ defined on the interval $[a,b]$, this formula is derived by introducing a linear combination of the derivative $f'$ computed at $n+1$ equally spaced points in $[a,b]$, together with the second-order derivatives $f''$ computed at the limit points $a$ and $b$.\\

\noindent  We proved that the corresponding remainder can be minimized for an {\em ad hoc} choice of the weights involved in this linear combination, and can be significantly 
smaller than the one obtained with the classical second order Taylor's formula.\sa
\noindent  Then, we considered two usual applications of this Taylor-like expansion: the interpolation error and the numerical quadrature formula. We showed that using this 
approach improves both the Lagrange $P_2$- interpolation error estimate and the error bound of the Simpson rule in numerical integration.\sa
\noindent For the interpolation error, we showed that the upper bound of the errors we obtained is almost two times more precise than the one obtained by the classical 
Taylor formula. For the numerical integration, dealing with functions with only $C^{3}$ regularity, the new quadrature error based on Taylor-like formula was found to 
be bounded $1.33$ times less than the best one derived for another numerical integration rule.\sa
\noindent Other applications could also be concerned by this new second order Taylor-like formula. For instance, we could consider to improve the approximation error involved in 
ODE's approximation where Taylor's formula is the main tool used to derive numerical schemes.\sa
\noindent \textbf{\underline{Homages}:} The authors want to warmly dedicate this research to pay homage to the memory of Professors Andr\'e Avez and G\'erard Tronel who largely promote the passion of research and teaching in mathematics of their students.

\begin{thebibliography}{8}
%
%
\bibitem{Abdulle} A. Abdulle, G. Garegnani, A probabilistic finite element method based on random meshes: A posteriori error estimators and Bayesian inverse problems, {\em Comput. Methods Appl. Mech. Eng.}, 113961, pp. 384 (2021).
%
\bibitem{AsCh2014} F. Assous, J. Chaskalovic, Indeterminate Constants in Numerical Approximations of PDE's: a Pilot Study Using Data Mining Techniques, \emph{J. Comput Appl. Math.},  270, pp. 462-470, (2014).
%
\bibitem{Atki88} K.E Atkinson, \emph{An Introduction to Numerical Analysis}, 2nd Edition, John Wiley \& Sons, 1988.
%
\bibitem{Barnett_Dragomir} N.S. Barnett, S.S. Dragomir, Applications of Ostrowski's version of the Grüss inequality for trapezoid type rules, {\em  Tamkang J. Math.,} 37(2), PP. 163--173 (2006).
%
\bibitem{BuFa11} {R.L. Burden, D. Faires}, \emph{Numerical Analysis}, 9th Edition, Brooks/Cole, Pacific Grove, 2011.
%
\bibitem{Cerone} P. Cerone, S.S. Dragomir, Trapezoidal-type rules from an inequalities point of view, in: G. Anastassiou (Ed.), Handbook of
Analytic-Computational Methods in Applied Mathematics, CRC Press, New York, pp. 65--134 (2000).
%
\bibitem{ArXiv_JCH} J. Chaskalovic, A probabilistic approach for solutions of determinist PDE's as well as their finite element approximations, {\em Axioms}, 10, pp. 349 (2021).
%
\bibitem{ChaskaPDE} J. Chaskalovic, {\em Mathematical and numerical methods for partial differential equations}, Springer Verlag, 2013.
%
\bibitem{CMAM2} J. Chaskalovic, F. Assous, A new mixed functional-probabilistic approach for finite element accuracy, \emph{Computational Methods in Applied Mathematics}, DOI: https://doi.org/10.1515/cmam-2019-0089 (2019).
%
\bibitem{ChAs20} J. Chaskalovic, F. Assous, Explicit k-dependence for $P_{k}$ finite elements in $W^{m,p}$ error estimates: application to probabilistic laws for accuracy analysis, {\em Applicable Analysis}, DOI: 10.1080/00036811.2019.1698727 (2020).
%
\bibitem{MMA2021} J. Chaskalovic, F. Assous, , Numerical validation of probabilistic laws to evaluate finite element error estimates, \emph{Mathematical Modelling and Analysis},  26 (4), pp. 684--694 (2021).
%
\bibitem{ChAs2023} J. Chaskalovic, F. Assous, A refined first-order expansion formula in $\R^n$: Application to interpolation and finite element error estimates, \emph{submitted}, (2023).
%
\bibitem{arXiv_First_Order} J. Chaskalovic, H. Jamshidipour, A new first order expansion formula with a reduced remainder, Axioms, 11, 562, (2022).
%
\bibitem{Chen01} X.L. Cheng, Improvement of some Ostrowski-Gr\"uss type inequalities, {\em Computers Math. Applic.}, 42, pp.10--114 (2001).
%
\bibitem{ChSu02} X.L. Cheng, J. Sun, A Note on the Perturbed Trapezoid Inequality, {\em J. Inequal. Pure and  Appli. Math.}, 3-2, article 29 (2002).
%
\bibitem{Crouzeix_Mignot} M. Crouzeix, A.L. Mignot, {\em Analyse numérique des équations différentielles}, Seconde édition, Collection mathématiques appliquées pour la maîtrise, Masson, 1992.
%
\bibitem{DrACxx} S.S. Dragomir, R.P. Agarwal, P. Cerone, On Simpson's Inequality and Applications, {\em J. Inequal. Appl.}, 5, pp. 533--579 (2000).
%
\bibitem{Dragomir_Sofo} S.S. Dragomir, P. Cerone and A. Sofo, Some remarks on the trapezoid rule in numerical integration, Indian J. Pure Appl. Math., 31(5) (2000), 475-494.
%
\bibitem{DrWa97} S.S. Dragomir, S. Wang, An inequality of Ostrowski-Gr\"uss type and its applications to the estimation of error bounds for some special means and for some numerical quadrature rules, {\em Computers Math. Applic.}, 33 (11), pp. 15--20 (1997).
%
\bibitem{Hennig} P. Hennig, M.A. Osborne, M. Girolami, Probabilistic numerics and uncertainty in computations. {\em Proc. R. Soc. A Math. Phys. Eng. Sci.},  20150142 , pp 471 (2015).
%
\bibitem{Liu02} Z. Lui, A Inequality of Simpson type, {\em Proc. R. Soc. A}, 461, pp. 2155--2158 (2005).
%
\bibitem{MaPU00} M Matic, J. Pecaric, N. Ujevic, Improvement and further generalization of inequalities of Ostrowski-Gr\"uss type, {\em Computers Math. Applic.}, 39 (3-4), pp. 161--175 (2000).
%
\bibitem{Oates} C.J. Oates, T.J. Sullivan, A modern retrospective on probabilistic numerics, {\em Stat. Comput.}, 29, pp.1335--1351 (2019).
%
\bibitem{RaTho82} P.A. Raviart et J.M. Thomas, Introduction \`a l'analyse num\'erique des \'equations aux d\'eriv\'ees partielles, Masson (1982).
%
\bibitem{SuMa03} E. S\"uli, D. Mayers, {\em An Introduction to Numerical Analysis}, Cambridge University Press (2003).
%
\bibitem{Taylor} B. Taylor, Methodus incrementorum directa and inversa, Innys: London, UK, Prop.VII, Th.III, p. 21 (1717).
%
\end{thebibliography}
\end{document}